\documentclass[a4 paper,11pt]{amsart}

\usepackage{amsmath} %
\usepackage{amssymb} %
\usepackage{amscd} %
\usepackage{amsthm} %
\usepackage{mathrsfs} %
\usepackage{euscript}
\usepackage{eufrak}
\usepackage{ulem}
\usepackage[alphabetic]{amsrefs}
\usepackage{ascmac} % 枠付き
\usepackage{comment}
\usepackage{bm}%太字

\usepackage[dvipdfmx]{graphicx} 

\usepackage{hyperref}
\usepackage{xcolor}
\definecolor{unbleu}{rgb}{0.03, 0.15, 0.4}

\hypersetup{
pdfborder = {0 0 0},
colorlinks,
linkcolor=unbleu,
citecolor=unbleu,
urlcolor=unbleu
}

\usepackage{fancyhdr}
 \newtheorem{theorem}{Theorem}[section]
 \newtheorem{lemma}[theorem]{Lemma}

\newtheorem{maintheorem}{Theorem} %\theoremstyle{definition}の下におくと斜体にならない!
  
\theoremstyle{definition}
\newtheorem{definition}[theorem]{Definition}

\def\ul#1{\underline{#1}}%  underline
\begin{document}

\title[]{Hausdorff dimension of the parameters for $(\alpha,\beta)$-transformations with the specification property}
%"Constraint Ergodic" -> "Constrained ergodic" 

\author[M. Oguchi]{Mai Oguchi}
\address{Department of Mathematics, Ochanomizu University, 2-1-1 Otsuka, Bunkyo-ku, Tokyo, 112-8610, Japan}
\email{}

\author[M. Shinoda]{Mao Shinoda}
\address{Department of Mathematics, Ochanomizu University, 2-1-1 Otsuka, Bunkyo-ku, Tokyo, 112-8610, Japan}
\email{shinoda.mao@ocha.ac.jp}

% \date{\today}
\subjclass[2010]{\textcolor{black}{Primary} 37E45, 37E05, 11K55}
%\textcolor{black}{; Secondary 37J51}}
\keywords{}

%%%%%%%%%%%%%%%%%%%%%%%%%%%%%%%%%%%%%%%%%%%%%%%%%%%%%%%%%
%%%%%%%%%%%%%%%%%%%%%%%%%%%%%%%%%%%%%%%%%%%%%%%%%%%%%%%%%
\begin{abstract}
In this paper we consider the specification property for $(\alpha,\beta)$-shifts.
When $\alpha=0$, Schmeling shows that the set of $\beta>1$ for which the $\beta$-shift has the specification property has the Lebesgue measure zero but has the full Hausdorff dimension\cite{Schmeling}.
So it is natural to ask what happens when $\alpha>0$.
Buzzi shows that for fixed $\alpha$ the set of $\beta >1$ for which the $(\alpha,\beta)$-shift has the specification properyt has Lebesgue measure zero.
Hence we consider the Hausdorff dimension of the parametar space of $(\alpha,\beta)$-shifts.
\end{abstract}

\maketitle
% %\newpage

% \tableofcontents
%%%%%%%%%%%%%%%%%%%%%%%%%%%%%%%%%%%%%%%%%%%%%%%%%%%%%%%%
\section{Introduction}

In this paper we consider the specification property for $(\alpha,\beta)$-shifts, which is the coding spaces of $(\alpha,\beta)$-transformations.
More precisely, for $0\leq \alpha<1$ and $\beta>1$ the $(\alpha,\beta)$-transformation $T_{\alpha,\beta}: [0,1)\rightarrow [0,1)$ is defined by
\begin{equation}
T_{\alpha,\beta}(x)=\beta x+\alpha-\lfloor\beta x+\alpha\rfloor%=\beta x+\alpha\;\;\;\mod1
\label{alpha_beta}
\end{equation}
 for all $x\in (0,1]$ where $\lfloor \xi \rfloor$ denotes the largest integer no more than $\xi$
 and the $(\alpha,\beta)$-shift $\Sigma_{\alpha,\beta}$ is a subshift consisting of all $(\alpha,\beta)$-expansions induced by the $(\alpha,\beta)$-transformation.
After the pioneering work of Renyi \cite{R57} and Parry \cite{Parry60,Parry64} dynamical properties of $T_{\alpha,\beta}$ and $\Sigma_{\alpha,\beta}$ have been investigated from view points of ergodic theory and number theory for over half a century \cite{Cara, Fall09, Rai11,Hof80,ChuYam23}.

The specification property, which is introduced by Bowen\cite{Bow71}, is very important since it allows the study of generic properties of invariant measures.
For example, it implies the uniqueness of equilibrium measure for any smooth potential with the additional hypothesis of expansiveness.

% It is natural to pay attention to the parameters $(\alpha,\beta)\in [0,1)\times (1,\infty)$ for which $\Sigma_{\alpha,\beta}$ has the specification property.
In the case of $\alpha=0$, the shift is known as the $\beta$-shift $\Sigma_\beta:=\Sigma_{0,\beta}$ and the parameter $\beta>1$ for which $\Sigma_\beta$ has the specification property is studied: the set of such parameters has the Lebesgue measure $0$ but has the full Hausdorff dimension \cite{Schmeling}.
For a fixed $\alpha\in [0,1)$ it is also known that $\beta>1$ for which $\Sigma_{\alpha,\beta}$ has the specification property has the Lebesgue measure $0$ \cite{Buzzi97}.
Hence it is natural to ask the Hausdorff dimension in the case of $\alpha\in (0,1)$.

It is well known that the $(\alpha,\beta)$-shift consists of the $(\alpha,\beta)$-expansions in the interval with regard to the lexicographical order $\preceq$ whose end points are $(\alpha,\beta)$-expansions of $0$ and $1$ (See for more details \S\ref{Preliminaries}).
Let $\mathbb{N}_0$ be the set of non negative integers.
Fixing the expansion of $0$, we show that the set of the parameters for which the $(\alpha,\beta)$-shift has  the specification property has the full Hausdorff dimension:
\begin{maintheorem}
\label{maintheorem}
    Let $\ul{u}\in \mathbb{N}_0^{\mathbb{N}_0}$ satisfy 
    $u_0=0$, $\ul{u}\preceq \sigma^n \ul{u}$ for all $n\geq 0$ and $\max_{i\in \mathbb{N}_0}u_i=:K<\infty$.
    Then we have 
    \begin{align*}
        \dim_H\{\beta>K+2: \Sigma_{\alpha(\beta),\beta}\ \mbox{has the specification property}\}=1
    \end{align*}
    where $\displaystyle \alpha(\beta)=(\beta-1)\sum_{j=0}^\infty \frac{u_j}{\beta^{j+1
}}$ and $\dim_H E$ is the Hausdorff dimension of $E$.
\end{maintheorem}
Our proof is based on the work by Hu et.al \cite{Hu}, which shows the set of $\beta>1$ for which the $(-\beta)$-shift has the specification property has the full Hausdorff dimension. 
The main tool is the well know theorem on the Hausdorff dimension of the attractor for an interacted function system. 

The rest of this paper organized as follows. 
In \S2, we establish our definitions and prepare several facts. 
In \S3 we prove our main theorem.

% many authors have extensively investigated $(\alpha,\beta)$-transformations and the associated $(\alpha,\beta)$-shifts.
% This paper investigates a dynamical system $(X,T)$ where $X$ is a compact metric space endowed with its Borel $\sigma$-alpgebra $\mathcal{B}$ and $T:X\rightarrow X$ is a measureble map.
% R. Bowein \cite{B} has introduced the specification property: a map $T$ has the specification property if
% for every $\varepsilon>0$ there exists a constant $D<\infty$ such that given $(x_1, \ell_1),\ldots, (x_r, \ell_r)\in X\times \mathbb{N}$, $0\leq r<\infty$, if we let $L(s)=\ell_1+\ldots+\ell_{s-1}+(s-1)D$ for $1\leq s \leq r+1$, then there exits $z\in X$ such that $f^{L(r+1)}(z)=z$ and $d(f^{L(s)+k}(z), f^{k}(x_s))<\varepsilon$ for $0\leq k<\ell_s$.
% The dynamical behaviors exhibited by $T_{\alpha,\beta}$ and $\Sigma_{\alpha,\beta}$ are of great significance in both ergodic theory and number theory.
%This map is related to the expansion of numbers %Fallerの言い方参照．

%%%%%%%%%%%%%%%%%%%%%%%%%%%%%%%%%%%%%%%%%%%%%%%%%%%%%%%%
\section{Preliminaries}\label{Preliminaries}
\subsection{Symbolic dynamics}
Let $A$ be a finite set and $A^{\mathbb{N}_0}$ the set of all one-sided infinite sequences on the alphabet $A$
with the standard metric $d(x,y)=2^{-t(x,y)}$, where $t(x,y)=\min\{k\in\mathbb{N}_0:x_k\not=y_k\}$.
Then the shift map $\sigma\colon A^{\mathbb{N}_0}\to A^{\mathbb{N}_0}$ defined by
$\sigma(x)_i=x_{i+1}$ is a continuous map. If a compact subset $X\subset A^{\mathbb{N}}$
satisfies $\sigma(X)\subset X$, then $X$ is called a {\it subshift}. 
For a finite sequence $w\in \bigcup_{n\ge 1}A^n$, denote by $w^{\infty}:=ww\cdots\in A^{\mathbb{N}}$
the periodic infinite sequence and by
$[w]:=\{x\in X:x_1\cdots x_{|w|}=w\}$ the \textit{cylinder set} for $w$, where $|w|$ denotes the length of $w$.
The \textit{language} $\mathcal{L}(X)$ is defined by
$$\mathcal{L}(X):=\left\{w\in\bigcup_{n\ge 1}A^n:[w]\not=\emptyset\right\}.$$
The topological entropy of $X$ is defined by
\begin{align*}
    h_{{\rm top}}(X)=\lim_{n\to \infty}\frac{1}{n}\log \#\mathcal{L}(X).
\end{align*}
% For a collection $\mathcal{D}\subset\mathcal{L}(X)$, we define the \textit{entropy} of $\mathcal{D}$ by
% $$h(\mathcal{D}):=\limsup_{n\rightarrow\infty}\frac{1}{n}\log\sharp\mathcal{D}_n,$$
% where $\mathcal{D}_n:=\{w\in\mathcal{D}:|w|=n\}$, and write $h_{\rm top}(X):=h(\mathcal{L}(X))$.
% The \textit{entropy of an invariant measure} $\mu$ is given by
% $$h(\mu):=\lim_{n\rightarrow\infty}\frac{1}{n}\sum_{w\in\mathcal{L}_n(X)}-\mu[w]\log \mu[w].$$
% Then well-known variational principle \cite[Theorem 8.6]{W2} implies that
% $$h_{\rm top}(X)=\sup\{h(\mu):\mu\text{ is invariant}\}.$$
% An invariant measure which attains this supremum is called a {\it measure of maximal entropy}.
% If such a measure exists and unique, then we say that $(X,\sigma)$ is {\it intrinsically ergodic}.

% As we said before, we use Climenhaga and Thompson's technique found in \cite{CT,CT2}
% to prove Theorem \ref{main1}.
% Before stating their results, 
% %we give definitions and notations which are appeared in them.
% %First, 
% we give definitions
%  of ``non-uniform" version of
% specification and Gibbs property.

\begin{definition}[Specification]
	Let $\mathcal{L}(X)$ be the language of a subshift $X$.
	We say that $X$ has {\it specification} with gap size $t\in \mathbb{N}$ if
	\begin{itemize}
		\item 
		for all $m\in \mathbb{N}$ and $w^1, \ldots, w^m\in \mathcal{L}(X)$, 
		there exist $v^1, \ldots, v^{m-1}\in \mathcal{L}(X)$ such that 
		$x:=w^1 v^1 w^2 v^2 \cdots v^{m-1} w^m \in \mathcal{L}(X)$
		and $|v^i|=t$ for all $i$ and
		
		\item
		the cylinder $[x]$ contains a periodic point of period exactly $|x|+t$.
		
	\end{itemize}
\end{definition}

 %%%%%%%%%%%%%%%%%%%%%%%%%%%%%%%%%%%%%%%%%%%%%%%%%%%%%%%%
\subsection{$(\alpha,\beta)$-shifts}

Let $\alpha\in[0,1)$ and $\beta>1$ and $T_{\alpha,\beta}$ be the $(\alpha,\beta)$-transformation defined by \eqref{alpha_beta}.
Let $k=\lceil\alpha+\beta\rceil$ and $A=\{0,1, \ldots, k-1\}$.
Set $I_0=\left[0, \frac{1-\alpha}{\beta}\right)$, $I_j=\left[\frac{j-\alpha}{\beta}, \frac{j+1-\alpha}{\beta}\right)$ for $1\leq j\leq k-2$ and $I_{k-1}=\left[\frac{k-1-\alpha}{\beta}\right)$.
 % \begin{align*}
 %     I_0=\left[0, \frac{1-\alpha}{\beta}\right)
 %     \quad\mbox{and}\quad
 %     I_{k-1}=\left[\frac{k-1-\alpha}{\beta}\right)
 % \end{align*}
 % and for $j=1, \ldots,k-2$ set
 % \begin{align*}
 %     I_j=\left[\frac{j-\alpha}{\beta}, \frac{j+1-\alpha}{\beta}\right).
 % \end{align*}
 Define the coding map $i_{\alpha,\beta}: [0,1)\rightarrow A^{\mathbb{N}_0}$ by $i_{\alpha,\beta}(x)_n=j$ iff $ T^n_{\alpha, \beta}(x)\in I_j$
for all $x\in [0,1)$ and $n\geq 0$.
%  \begin{align*}
%     i_{\alpha,\beta}(x)_n=j \Leftrightarrow T^n_{\alpha, \beta}(x)\in I_j
% \end{align*}
Then we have
\begin{align*}
    x=\frac{i_{\alpha,\beta}(x)_0-\alpha}{\beta}+\frac{i_{\alpha,\beta}(x)_1-\alpha}{\beta^2}+\cdots=\sum_{n=0}^\infty \frac{i_{\alpha,\beta}(x)_n-\alpha}{\beta^n}.
\end{align*}
We call $i_{\alpha,\beta}(x)$ a $(\alpha,\beta)$-expansion of $x$.
Now we define
\begin{align*}
    \Sigma_{\alpha,\beta}:={\rm cl}(\{i_{\alpha,\beta}(x)\in A^{\mathbb{N}_0}: x\in (0,1]\})
\end{align*}
where ${\rm cl}(X)$ is the clusure of $X$ in $A^{\mathbb{N}_0}$.
Then it is easy to see that $\Sigma_{\alpha,\beta}$ is a subshift of $A^{\mathbb{N}}$.
We call $(\Sigma_{\alpha,\beta}, \sigma)$ the $(\alpha,\beta)$-shift.
It is well known that the $(\alpha,\beta)$-shift are characterized by the lexicographical order and the expansions of the end points: 
\begin{align}
        \Sigma_{\alpha,\beta}=\{\ul{x}\in A^{\mathbb{N}_0}: \ul{u}^{\alpha,\beta}\preceq \sigma^n \ul{x}\preceq \ul{v}^{\alpha,\beta}\ \mbox{for all}\ n\geq 0\}
        \label{interval}
    \end{align}
where $\ul{u}^{\alpha,\beta}=i_{\alpha,\beta}(0)$ and $\displaystyle \ul{v}^{\alpha,\beta}=\lim_{x \uparrow 1}i_{\alpha,\beta}(x)$.
 Moreover we have
\begin{align*}
    \ul{u}^{\alpha,\beta}\preceq\sigma^n \ul{u}^{\alpha,\beta}\preceq \ul{v}^{\alpha,\beta}
    \quad \mbox{and}\quad 
    \ul{u}^{\alpha,\beta}\preceq\sigma^n \ul{v}^{\alpha,\beta}\preceq \ul{v}^{\alpha,\beta}
\end{align*}
for all $n\geq 0$.
See \cite{Faller} and references therein  for more details.

% \begin{proposition}\cite{Faller} %主張を発見する
% \label{interval}
%     Let $\Sigma_{\alpha,\beta}$ be a $(\alpha,\beta)$-shfit.
%     Then we have
%     \begin{align*}
%         \Sigma_{\alpha,\beta}=\{\ul{x}\in A^{\mathbb{N}_0}: \ul{u}^{\alpha,\beta}\preceq \sigma^n \ul{x}\preceq \ul{v}^{\alpha,\beta}\ \mbox{for all}\ n\geq 0\}
%     \end{align*}
%     where $\ul{u}^{\alpha,\beta}=i_{\alpha,\beta}(0)$ and $\displaystyle \ul{v}^{\alpha,\beta}=\lim_{x \uparrow 1}i_{\alpha,\beta}(x)$.

%     Moreover we have
% \begin{align*}
%     \ul{u}^{\alpha,\beta}\preceq\sigma^n \ul{u}^{\alpha,\beta}\preceq \ul{v}^{\alpha,\beta}
%     \quad \mbox{and}\quad 
%     \ul{u}^{\alpha,\beta}\preceq\sigma^n \ul{v}^{\alpha,\beta}\preceq \ul{v}^{\alpha,\beta}
% \end{align*}
% for all $n\geq 0$.
% \end{proposition}

For $(\alpha,\beta)$-shifts the specification property is characterized by the expansions $\ul{u}^{\alpha,\beta}$ and $\ul{v}^{\alpha,\beta}$.
For $\ul{u}^{\alpha,\beta},\ul{v}^{\alpha,\beta}$ define
\begin{align*}
D(\ul{u}^{\alpha,\beta})=\{n\in\mathbb{N}\;:\;(v^{\alpha,\beta}_0,\dots,v^{\alpha,\beta}_n)=(u^{\alpha,\beta}_j,\dots,u^{\alpha,\beta}_{j+n})\;\text{for some}\;j\in\mathbb{N}\}\\
D(\ul{v}^{\alpha,\beta})=\{n\in\mathbb{N}\;:\;(u^{\alpha,\beta}_0,\dots,u^{\alpha,\beta}_n)=(v^{\alpha,\beta}_j,\dots,v^{\alpha,\beta}_{j+n})\;\text{for some}\;j\in\mathbb{N}\}.
\end{align*}
%Then we have the following characterization of the specification property.
\begin{theorem}\cite[Theorem 1.5]{Cara}
\label{specification}
    Let $0\leq \alpha<1$ and $\beta>2$.
    Then $\Sigma_{\alpha,\beta}$ has the specification property if and only if both $D(\ul{u}^{\alpha,\beta})$ and $D(\ul{v}^{\alpha,\beta})$ are bounded.
\end{theorem}

% %%%%%%%%%%%%%%%%%%%%%%%%%%%%%%%%%%%%%%%%%%%%%%%%%%%%%%%%
% \subsection{The inverse problem for $(\alpha,\beta)$-shifts}

% Since the $(\alpha,\beta)$-shifts defined by the two sequences $u^{\alpha,\beta}, v^{\alpha,\beta}$ as in \eqref{interval}, it is natural to ask the inverse problem:
% given two sequences $\ul{u}$ and $\ul{v}$ satisfying
% \begin{align}
%     \ul{u}\preceq\sigma^n \ul{u}\preceq \ul{v}
%     \quad \mbox{and}\quad 
%     \ul{u}\preceq\sigma^n \ul{v}\preceq \ul{v}
%     \quad ^\forall n\geq 0.
%     \label{inequality}
% \end{align}
% Can we find $\alpha\in [0,1)$ and $\beta \in (1,\infty)$ s.t. $\ul{u}=\ul{u}^{\alpha,\beta}$ and $\ul{v}=\ul{v}^{\alpha,\beta}$?
% The answer is given in \cite{Faller}.
% Let $k\geq 3$ and let $\ul{u}\in \{0,1, \ldots, k-1\}^{\mathbb{N}_0}$ and $\ul{v}\in \{0,1, \ldots, k-1\}^{\mathbb{N}_0}$ s.t. $u_0=0, v_0=k-1$ and \eqref{inequality} holds.
% Set $\hat{\ul{u}}=\max\{\sigma^n\ul{u}: n\geq 0\}$ and 
% \begin{align*}
%     \Sigma(\ul{u},\hat{\ul{u}})=\{\ul{x}\in \{0,1, \ldots, k-1\}^{\mathbb{N}_0}: \ul{u}\preceq \sigma^n \ul{x}\preceq \hat{\ul{u}}\ \mbox{for all}\ n\geq 0\}.
% \end{align*}
% By Theorem 6.20 in \cite{Faller} there exist $\alpha\in [0,1)$ and $\beta >1$ such that $\ul{u}=\ul{u}^{\alpha,\beta}$ and $\ul{v}=\ul{v}^{\alpha,\beta}$, 
% if $k-2> h_{{\rm top}}(\Sigma(\ul{u}, \hat{\ul{u}}))$.
% We use this fact later.

%%%%%%%%%%%%%%%%%%%%%%%%%%%%%%%%%%%%%%%%%%%%%%%%%%%%%%%%
\section{Proof of the main theorem}
Our proof of Theorem \ref{maintheorem} is based on the idea by Hu et.al \cite{Hu}, which proves that the set of $\beta>1$ for which the $(-\beta)$-shift has the specification property has full Hausdorff dimension.
 
Let $\ul{u}\in \mathbb{N}^{\mathbb{N}_0}$ satisfy the following:
\begin{description}
    \item[U1] $u_0=0$
    \item[U2] $\ul{u}\prec \sigma^n \ul{u}$ for every $n\geq 1$ and
    \item[U3] $\displaystyle \max_{i\in \mathbb{N}_0} u_i=:K<\infty$
\end{description}
Since the case $\ul{u}=0^\infty$ corresponding to the case $\alpha=0$ is proved by Schmeling\cite{Schmeling},  we may assume $\ul{u}\neq 0^\infty$.
In the lest of this paper we assume $ N\geq K+3$.

\begin{lemma}
\label{onetoone}
    Let $\ul{v}\in \{0,1, \ldots, N-1\}^{\mathbb{N}_0}$ such that
    $v_0=N-1$ and $1\leq v_i\leq N-2$ for all $i\geq 1$.
    Then there exist $(\alpha,\beta)\in (0,1)\times (1, \infty)$ satisfying
    \begin{align*}
        0=\sum_{i=0}^\infty \frac{u_i-\alpha}{\beta^{i+1}}
        \quad \mbox{and}\quad
            1=\sum_{i=0}^\infty \frac{v_i-\alpha}{\beta^{i+1}}.
    \end{align*}
    Moreover we have $\ul{u}^{\alpha,\beta}=\ul{u}$ and $\ul{v}^{\alpha,\beta}=\ul{v}$.

\begin{proof}
   Since $u_0=0$, $v_0=N-1$ and $1\leq v_i\leq N-2$ for all $i\geq 1$, it is easy to see that
\begin{align}
    \ul{u}\preceq\sigma^n \ul{u}\preceq \ul{v}
    \quad \mbox{and}\quad 
    \ul{u}\preceq\sigma^n \ul{v}\preceq \ul{v}
    \quad ^\forall n\geq 0
    \label{inequality}
\end{align}
holds.
Set $\hat{\ul{u}}=\max\{\sigma^n\ul{u}: n\geq 0\}$ and 
\begin{align*}
    \Sigma(\ul{u},\hat{\ul{u}})=\{\ul{x}\in \{0,1, \ldots, N-1\}^{\mathbb{N}_0}: \ul{u}\preceq \sigma^n \ul{x}\preceq \hat{\ul{u}}\ \mbox{for all}\ n\geq 0\}.
\end{align*}
By U3 we have $\Sigma(\ul{u}, \hat{\ul{u}})\subset \{0,1, \ldots, K\}^{\mathbb{N}_0}$ and $h_{top}(\Sigma(\ul{u},\hat{\ul{u}}))\leq \log (K+1)$.
By the choice of $N$, we have $\log v_0=\log K+2 >\log (K+1)=h_{{\rm top}}(\Sigma(\ul{u}, \hat{\ul{u}}))$.
Then by Theorem 6.20 in \cite{Faller} there exist $\alpha\in (0,1)$ and $\beta >1$ such that $\ul{u}=\ul{u}^{\alpha,\beta}$ and $\ul{v}=\ul{v}^{\alpha,\beta}$.
\end{proof}
\end{lemma}
% This lemma shows that there is a one-to-one correspondence between the pair $(\alpha,\beta)$ and $(\ul{u},\ul{v})$.
Let 
\begin{align*}
    E_N=\{\beta>N-2:  v_0^{\alpha(\beta), \beta}=N-1, 1\leq v_i^{\alpha(\beta),\beta}\leq N-2 \ \mbox{for all}\ i\geq 1\ \}
\end{align*}
where $\displaystyle \alpha(\beta)=(\beta-1)\sum_{j=0}^\infty \frac{u_j}{\beta^{j+1}}$.
For $\beta \in E_N$ we have $D(\ul{u}^{\alpha(\beta), \beta})=D(\ul{v}^{\alpha(\beta),\beta})=\emptyset$
since $v_0^{\alpha(\beta),\beta}=N-1\neq u_i$ for all $i\in \mathbb{N}_0$ and $u_0^{\alpha(\beta),\beta}=0\neq v_i$ for all $i\in \mathbb{N}_0$.
By Theorem\ref{specification} the shift $\Sigma_{\alpha(\beta),\beta}$ satisfies the specification property.

For $\beta>2$ we have
\begin{align*}
    \frac{d \alpha(\beta)}{d\beta}= \sum_{j=1}^\infty \left(-j+\frac{1}{\beta}(j+1)\right)<0,
\end{align*}
which implies $\alpha(\beta)$ is strictly decreasing.
By Lemma \ref{onetoone} for $\beta \in E_N$ we have $\beta\leq N$ and
\begin{align*}
    \alpha(\beta)\leq  \alpha(N-2)\leq \frac{K}{N-2}=:\alpha_N<1.
\end{align*}
Hence we have $E_N\subset (N-1-\alpha_N,N]$.

Define a map $\varphi: E_N\rightarrow [0,1]$ by 
\begin{align*}
    \varphi(\beta)=\sum_{k=1}^\infty \frac{v_k^{\alpha,\beta}}{N^k}.
\end{align*}
For $i=1, \ldots, N-2$ set $\displaystyle \psi_i(x)=\frac{x+i}{N}$.
Then by Lemma \ref{onetoone} it is easy to see $\varphi(E_N)$ satisfies
\begin{align*}
    \varphi(E_N)=\bigcup_{i=1}^{N-2} \psi_i(\varphi(E_N))
\end{align*}
and  it is the attractor of the interacted function system $\Psi=\{\psi_i\}_{i=1}^{N-2}$.
By \cite[Theorem 9.3]{Fal90} $\dim_H \varphi(E_N)=s$, where $s$ is given by $\displaystyle \sum_{i=1}^{N-2}\left( \frac{1}{N}\right)^s=1$.
Hence we have
\begin{align*}
    \dim_H \varphi(E_N)=\frac{\log(N-3)}{\log N}.
\end{align*}

\begin{lemma}
\label{Lip}
    There exists $C_N>0$ such that for every $\beta, \beta'\in E_N$ we have
    \begin{align*}
        |\varphi(\beta)-\varphi(\beta')|\leq C_N |\beta-\beta'|.
    \end{align*}
\end{lemma}

We will prove this lemma later.
By using this lemma and \cite{Fal90}[Proposition 3.3] we have
\begin{align*}
    \dim_H E_N\geq \dim_H \varphi(E_N)=\frac{\log(N-3)}{\log N}.
\end{align*}
Since this holds for ever $N\geq K+3$, Theorem \ref{maintheorem} is proved.

In the rest of this paper, we prove Lemma \ref{Lip}.
Let $\beta, \beta'\in E_N$ and $n\geq 1$ satisfy
\begin{align*}
    v_i^{\alpha(\beta),\beta}=v_i^{\alpha(\beta'),\beta'}
    \quad(0\leq i\leq n),\quad
    v_{n+1}^{\alpha(\beta), \beta}\neq v_{n+1}^{\alpha(\beta'),\beta'}.
\end{align*}
By the definition of $\varphi$ it is easy to see that
\begin{align}
    |\varphi(\beta)-\varphi(\beta')|\leq N^{-n}.
    \label{upper}
\end{align}

Without loss of generality we may assume $\beta<\beta'$.
For $0\leq i\leq n$ set $v_i:=v_i^{\alpha(\beta),\beta}=v_i^{\alpha(\beta'),\beta'}$.
Define polynomials by $  P_0(x)=x+\alpha(x)-v_0$ and 
\begin{align*}
    P_k(x)=xP_{k-1}(x)-\alpha(x)-v_k
\end{align*}
for $1\leq k \leq n$.
Then for $0\leq k\leq n$ we have $T^k_{\alpha(\beta),\beta}(1)=P_k(\beta)$ and $T^k_{\alpha(\beta'),\beta'}(1)=P_k(\beta')$.
Since
\begin{align*}
    P'_0(x)=1+\alpha'(x),
    \quad 
    P'_k(x)=P_{k-1}(x)+xP'_{k-1}(x)+\alpha'(x)
    \ (1\leq k\leq n),
\end{align*}
 we have
 \begin{align}
     P'_n(x)=x^{n-1}\left((1+\alpha'(x))x+\sum_{j=0}^{n-1}\frac{P_j(x)+\alpha'(x)}{x^j}\right).
     \label{derivative}
 \end{align}
By the mean value theorem there exists $\xi\in (\beta, \beta')\subset (N-1-\alpha_N,N]$ such that
\begin{align*}
    \left|T_{\alpha(\beta),\beta}^n(\beta)-T_{\alpha(\beta'),\beta'}^n(1)\right|
    &=\left|P_n(\beta)-P_n(\beta')\right|\\
    &=|\beta-\beta'|P'_n(\xi).
\end{align*}
By \eqref{derivative} and $P_j\in [0,1]$ for $0\leq j\leq n-1$ we have
\begin{align*}
    \left|T_{\alpha(\beta),\beta}^n(\beta)-T_{\alpha(\beta'),\beta'}^n(1)\right|
    &=|\beta-\beta'|\left((1+\alpha'(\xi))\xi^n+\sum_{j=0}^{n-1}\frac{P_j(\xi)+\alpha'(\xi)}{\xi^j}\right)\\
    &\leq |\beta-\beta'|\left((1+\alpha'(\xi))\xi^n+\left(1+\alpha'\left(\zeta\right)\right)\frac{1-\frac{1}{\zeta^n}}{1-\frac{1}{\zeta}}\right)\\
    &\leq |\beta-\beta'|\frac{\zeta^{2n}+\zeta^n-|\alpha'(\zeta)|(-\zeta^{2n-1}-1)}{\zeta^n-\zeta^{n-1}}\\
    &\leq|\beta-\beta'|\frac{2\zeta^{2n}+2|\alpha'(\zeta)|\zeta^{2n}}{\zeta^n-\zeta^{n-1}}\\
%    &=|\beta-\beta'|\frac{2\zeta^{n+1}+2|\alpha'(\zeta)|\zeta^{n+1}}{\zeta-1}\\
    &\leq|\beta-\beta'|\zeta^{n+1}(2+2|\alpha'(\zeta)|)\\
    &<|\beta-\beta'|N^{n+1}(2+2|\alpha'(\zeta)|).
\end{align*}
Therefore we have
\begin{align}
    |\beta-\beta'|> \left|T_{\alpha(\beta),\beta}^n(\beta)-T_{\alpha(\beta'),\beta'}^n(1)\right|\frac{1}{N^{n+1}(2+2|\alpha'(\zeta)|)}.
    \label{lower}
\end{align}

For $x,y\in (N-1-\alpha_N,N]$ define
\begin{align*}
\Delta_N(x,y)=\min\left\{\left(\frac{2-\alpha(x)}{x}+\displaystyle\sum^{\infty}_{k=2}\frac{1-\alpha(x)}{x^{k}}\right)-\left(\frac{1-\alpha(y)}{y}+\displaystyle\sum^{\infty}_{k=2}\frac{N-2-\alpha(y)}{y^{k}}\right),\right.\\
\left.\left(\frac{N-2-\alpha(y)}{y}+\displaystyle\sum^{\infty}_{k=2}\frac{1-\alpha(y)}{y^{k}}\right)-\left(\frac{N-3-\alpha(x)}{x}+\displaystyle\sum^{\infty}_{k=2}\frac{N-2-\alpha(x)}{x^{k}}\right)\right\}.
\end{align*}
Then we have
\begin{equation*}
\Delta_N(x,x)=\frac{1}{x}-\displaystyle\sum^{\infty}_{k=2}\frac{N-3}{x^k}=\frac{1}{x}\left(1-\frac{N-3}{x-1}\right)>\frac{1}{N}\left(1-\frac{N-3}{N-2-\alpha_N}\right)=:C_N>0.
\end{equation*}
Since $\Delta_N$ is continuous, there exists $\delta(N)>0$ such that for $x,y>2$ with $0\leq |x-y|\leq \delta(N)$ implies $\Delta(x,y)\geq C(N)>0$.

By the definition of the $(\alpha,\beta)$-transformation we have
\begin{equation*}
T^n_{\alpha(\beta),\beta}(1)=\displaystyle\sum^{\infty}_{i=0}\frac{v_{n+i}^{\alpha(\beta),\beta}-\alpha(\beta)}{\beta^{i+1}}\;\;,\;\;T^n_{\alpha(\beta'),\beta'}(1)=\displaystyle\sum^{\infty}_{i=0}\frac{v_{n+i}^{\alpha(\beta'),\beta'}-\alpha(\beta')}{\beta'^{i+1}}
\end{equation*}
%Hence by letting $0\leq \beta-\beta'\leq \delta(N)$ if necessary, we have
and by the definition of $\Delta_N$ we have
\begin{equation}\label{D_relation}
|T^n_{\alpha(\beta),\beta}(1)-T^n_{\alpha(\beta'),\beta'}(1)|\geq\Delta_N(\beta,\beta').
\end{equation}
By \eqref{upper},\eqref{lower} and \eqref{D_relation} we have
\begin{align*}
    |\beta-\beta'|\geq C_1(N) N^{-n} \geq C_1(N) |\varphi(\beta)-\varphi(\beta')|
\end{align*}
where $C_1(N)=C_N (N(2+2|\alpha'(\zeta)|))^{-1}$, which completes the proof.

%----------------------------------------------
% **********************************************************
% Appendices
% **********************************************************

% \setcounter{equation}{0}
% \renewcommand{\theequation}{\Alph{section}.\arabic{equation}}

% \appendix

% %\begin{comment}
% \section{On positive entropy}\label{AppendixA}

%%%%%%%%%%%%%%%%%%%%%%%%%%%%%%%%%%%%%%%%%%%%%%%%
\vspace*{33pt}

\noindent
\textbf{Acknowledgement.}~ 
We thank Kanji Inui for communicating us results before publication.
The second author was partially supported by JSPS KAKENHI Grant Number 21K13816.

\vspace{11pt }
\noindent
\textbf{Data Availability.}~
Data sharing not applicable to this article as no datasets were generated or analyzed during the current study.

%%%%%%%%%%%%%%%%%%%%%%%%%%%%%%%%%%%%%%%%%%%%%%%%%%%%%%%%
% \begin{thebibliography}{}
% \bibitem[Bou01]{Bou01} LA CONDITION DE WALTERS
% \bibitem[Cha11]{Cha11} Zero-temperature limit of one-dimensional Gibbs states via renormalization: the case of locally constant potentials. 
% \bibitem[Jen01]{Jen01} Rotation, entropy, and equilibrium states
% \bibitem[Jen06]{Jen06} Every ergodic measure is uniquely maximizing
% \bibitem[GarLop2018]{GarLop2018} Functions for relative maximization
% \bibitem[DGS]{DGS} Ergodic theory on compact spaces
% \bibitem[Mor2010]{Mor2010} Ergodic optimization for generic continuous functions
% \bibitem[Lep14]{Lep14} Zero temperature and selection of maximizing measure, http://www.math.univ-brest.fr/perso/renaud.leplaideur/KATHMANDU.pdf
% \end{thebibliography}

\bibliographystyle{alpha}
\bibliography{alpha_beta_Hausdorff.bib}

\end{document}